\documentclass[nomath,noamsfonts]{amsart}
\newcommand{\rom}{\mathrm}
\def\N{{\rom I\kern-.1567em N}}
\def\R{{\rom I\kern-.1567em R}}
\def\C{{\rom C\kern-6.5pt  
          \vrule height 7.7pt width 0.4pt depth -0.5pt \phantom {.}}}             
\newtheorem{satz}{Theorem}
\newtheorem{lemma}[satz]{Lemma}
\newtheorem{proposition}[satz]{Proposition}
\newtheorem{coro}[satz]{Corollary}

\newcommand{\falle}{\;\;\forall}
\newcommand{\gen}{\rightarrow}
\newcommand{\nach}[1]{\stackrel{#1}\rightarrow}
\newcommand{\Norm}[1]{\Bigl\|#1\Bigr\|}
\newcommand{\norm}[1]{\|#1\|}

\newcommand{\betr}[1]{| #1  |}
\newcommand{\eing}[1]{|_{#1}}
\newcommand{\ebew}{\hfill{\rule{1.2ex}{1.2ex}}}
\newcommand{\bgl}{\begin{eqnarray}}
\newcommand{\bglst}{\begin{eqnarray*}}
\newcommand{\egl}{\end{eqnarray}}
\newcommand{\eglst}{\end{eqnarray*}}
\newcommand{\Pel}{Pe\l\-czy\'ns\-ki}
\newcommand{\leins}{l^{1}}
\newcommand{\Ref}[1]{(\ref{#1})}
\newcommand{\id}{{{\rom i}\rom{d}}}
\newcommand{\genwst}{\nach{w^{*}}}
\newcommand{\wst}{$w^{*}$}
\newcommand{\w}{$w$}
\newcommand{\mdE}{|\;}
\begin{document}
\title{A note on asymptotically isometric copies of $\leins$
and $c_0$}
\author{Hermann Pfitzner}
\address{Universit\'e d'Orl\'eans\\
BP 6759\\
F-45067 Orl\'eans Cedex 2\\France}
\email{pfitzner@labomath.univ-orleans.fr}
\begin{abstract}
Nonreflexive Banach spaces that are complemented in their bidual by
an L-projection - like preduals of von Neumann algebras
or the Hardy space $H^1$ - contain, roughly speaking,
many copies of $\leins$ which are very close to isometric copies.
Such $\leins$-copies are known to fail the fixed point property.
Similar dual results hold for $c_0$.
\end{abstract}
\maketitle
In \cite{DJLT} it is shown that an isomorphic $\leins$-copy
does not necessarily contain asymptotically isometric $\leins$-copies
although by James' classical distortion theorem it always contains
almost isomorphic $\leins$-copies. (For definitions see below.)
Within this context and the context of the fixed point property
Dowling and Lennard \cite{DowLen}
show that the presence of an asymptotic $\leins$-copy
makes a Banach space fail the fixed point property.
Then they prove that every
nonreflexive subspace of $L^1[0,1]$ fails the fixed point property
by observing that
the proof of a theorem of Kadec and \Pel\ \cite[Th.\ 6]{KadPel}
yields an  asymptotic $\leins$-copy inside such supspaces of $L^1[0,1]$.
Alspach's example \cite{Alsp} may be considered as an early forerunner
of these results.

In the present note
we modify a construction of Godefroy in order to
show that every nonreflexive subspace of any L-embedded Banach space
contains an asymptotic $\leins$-copy and thus, in particular, fails
the fixed point property.
Analogous results hold for $c_0$ and M-embedded spaces.\medskip\\
Let $(x_n)$ be
a sequence of nonzero elements in a Banach space $X$.\\
We say that
\begin{em}$(x_n)$ spans $\leins$ $r$-isomorphically\end{em}
or just \begin{em}isomorphically\end{em}
if there exists $r>0$ (trivially $r\leq 1$) such that
$r\Bigl(\sum_{n=1}^{\infty}\betr{\alpha_n}\Bigr)\leq
\Norm{\sum_{n=1}^{\infty} \alpha_n {x_n}}   \leq
        \sum_{n=1}^{\infty}\betr{\alpha_n}$
for all scalars $\alpha_n$.
We say that
\begin{em}$(x_n)$ spans  $\leins$ almost isometrically\end{em}
if it is such that there exists a sequence
$(\delta_m)$ in $[0,1[$ tending to $0$
such that
$(1-\delta_m)\Bigl(\sum_{n=m}^{\infty}\betr{\alpha_n}\Bigr)\leq
\Norm{\sum_{n=m}^{\infty} \alpha_n {x_n}}   \leq
         \Bigl(\sum_{n=m}^{\infty}\betr{\alpha_n}\Bigr)$
for all $m\in\N$.
Trivially the property of spanning $\leins$ almost isometrically
passes to subsequences.
Analogously, we say that \begin{em}$(x_n)$ spans $c_0$
almost isometrically\end{em}
if  there exists a sequence
$(\delta_m)$ as above such that
$(1-\delta_m)\sup_{m\leq n\leq m'}\betr{\alpha_n}
\leq
\Norm{\sum_{n=m}^{m'} \alpha_n {x_n}}
\leq
(1+\delta_m)\sup_{m\leq n\leq m'}\betr{\alpha_n}$
for all $m\leq m'$.\medskip

Recall that James distortion theorem \cite{Jam64}
(or \cite{Die-Seq}, \cite{LiTz1})
for $\leins$ and $c_0$
says that every isomorphic
copy of $\leins$ (of $c_0$) contains an almost isometric copy
of $\leins$ (of $c_0$).
\medskip

Following \cite{DJLT, DowLen, DowLenTur}
we say that {\em $(x_n)$ spans $\leins$
asymptotically isometrically} (or just that {\em $(x_n)$ spans $\leins$
asymptotically}) if there is a sequence $(\delta_n)$ in $[0,1[$
tending to $0$ such that
$\sum_{n=1}^{\infty}(1-\delta_n)\betr{\alpha_n}\leq
\Norm{\sum_{n=1}^{\infty} \alpha_n {x_n}}   \leq
         \sum_{n=1}^{\infty}\betr{\alpha_n}$
for all scalars $\alpha_n$.
A sequence $(x_n)$ is said to {\em span $c_0$
asymptotically isometrically} (or just {\em to span $c_0$
asymptotically})
if there is a sequence $(\delta_n)$ as above such that
$\sup_{n\leq m}(1-\delta_n)\betr{\alpha_n}
\leq \Norm{\sum_{n=1}^{m}\alpha_n {x_n}}
\leq \sup_{n\leq m}(1+\delta_n)\betr{\alpha_n}$
for all $m\in\N$.
(Note that in the definition of asymptotic $c_0$'s and $\leins$'s
in \cite{DJLT} the sequence $(\delta_n)$ is supposed to
be decreasing but that this difference is not essential.)
Finally we say that a Banach space is isomorphic
(respectively almost isometric, respectively asymptotically isometric)
to $\leins$ (to $c_0$) if it has a basis with the corresponding property.
Clearly a sequence spanning $\leins$ asymptotically spans
$\leins$ almost isometrically.
The main result of \cite{DJLT} states that the converse does not
hold because there are almost isometric copies of $\leins$
which do not contain $\leins$ asymptotically.
Analogously, it was proved in \cite[Th.\ 3]{DJLT} that there exist
almost isometric copies of $c_0$ which do not
contain asymptotic $c_0$-copies.

Some notation:
The results are stated for complex scalars but hold also for
real Banach spaces.
Operator means linear bounded map.
As usual, we consider a Banach
space as a subspace of its bidual omitting the canonical embedding.
$[x_n]$ denotes the complete linear hull of $(x_n)$, $(e_n)$ is
the standard basis of $\leins$.
Basic properties and definitions of Banach space theory
can be found in \cite{Die-Seq} or in \cite{LiTz1, LiTz2}.

A Banach space $X$ is said to have the {\em fixed point property} if any contractive
(not necessarily linear) map $f:C\gen C$ on any non-empty
closed bounded convex subset $C\subset X$
has a fixed point where contractive means that
$\norm{f(x)-f(y)}<\norm{x-y}$ for all $x,y\in X$.

Let $Y$ be a subspace of a Banach space $X$ and $P$ be a projection
on
$X$.
$P$ is called an \begin{em}L-projection\end{em}
provided
$\norm{x}=\norm{Px}+\norm{(\id_{X}-P)x}$ for all $x\in X$.
A subspace $Y\subset X$ is called an \begin{em}M-ideal in $X$\end{em}
if its annihilator $Y^{\bot}$ in $X'$ is the
range of an L-projection on $X'$.
$Y$ is called an \begin{em}L-summand in $X$\end{em} if it is the
range of an
L-projection on $X$.
In the special case where $X=Y''$ and where $Y$ is an M-summand
(respectively an L-summand) in $Y''$  we say
that $Y$ is \begin{em}M-embedded\end{em} (respectively
\begin{em}L-embedded\end{em}).
As examples we only mention that preduals of von Neumann algebras,
in particular $\leins$ and $L^1$-spaces,
furthermore the Hardy space $H_0^1$ and the dual
of the disc algebra are L-embedded.
The sequence space $c_0$,
the space of compact operators on a Hilbert space, and
the quotient $C/A$ of the
continuous functions on the unit circle by the disc algebra $A$
are examples among M-embedded spaces.
The dual of an M-embedded space is L-embedded;
the converse is false \cite[III.1.3]{HWW}.
Throughout this note, if $X$ is an L-embedded Banach space
we will write $X_s$ for
the complement of (the canonical embedding of) $X$ in $X''$
that is $X''=X\oplus_1 X_s$.
In this case $P$ (or $P_X$ to avoid confusion)
will denote the L-projection from $X''$ onto $X$.
There is a useful criterion for L-embeddedness of
subsapces of L-embedded spaces due to
Li (\cite{Li-Ox} or \cite[Th.\ IV.1.2]{HWW}):
\begin{em}
A closed subspace $Y$ of an L-embedded Banach space $X$
is L-embedded if and only if $P Y^{\bot\bot}=Y$.
\end{em}
In this case if $Y$ is L-embedded and if one identifies
$Y''=Y\oplus_1 Y_s$ and $Y^{\bot\bot}\subset X''$
then $Y_s=Y^{\bot\bot}\cap X_s$.
Since biduals of L-embedded spaces are quite "big" and therefore
difficult to handle we mention in passing that a theorem of
Buhvalov-Lozanovskii (\cite{Buh-Lo}, \cite[IV.3.4]{HWW})
provides a characterisation of
L-embeddedness of subspaces of $L^1[0,1]$
only in terms of the spaces themselves: Subspaces of $L^1[0,1]$
are L-embedded if and only if their unit balls are closed with respect
to the measure topology \cite[IV.3.5]{HWW}.\\
The standard reference for M- and L-embedded spaces is the monograph
\cite{HWW}.\medskip

There are few stability results for L-embeddedness: Neither subspaces
nor quotients inherit this property
\cite[IV.1]{HWW}.
Therefore the following lemma reveals a nice exception.
In particular it underlines the idea that $\leins$-copies are
the building blocks of L-embedded spaces.
\begin{lemma}\label{prop_LLasy}
Almost isometric copies of $\leins$ which are subspaces of
L-embedded Banach spaces are L-embedded.\end{lemma}
{\em Proof}:
Consider first the L-embedded subspaces $U_m=[e_n]_{n\geq m}$
of $\leins$, $m\in\N$. We have
$U_m^{\bot\bot}=U_m \oplus_1 (c_0 \cap U_m^{\bot\bot})$.
An element $\mu\in(\leins)_s=c_0^{\bot}$ belongs to each
$U_m^{\bot\bot}$.
[Denote by $\rho_m$ the projection
$(\alpha_n)\mapsto (0,\ldots,0,\alpha_{m+1},\alpha_{m+2})$ on $\leins$.
Then $\mu$ annihilates ${\mbox{ker}\,\rho'_m}$
because $\mbox{ker}\,\rho'_m\subset c_0$.
Thus $\mu\in\mbox{ran}(\rho'')$ and
$\mu\in\overline{[e_n]}^{w^*}_{n\geq m}\cap c_0^{\bot}$.]

Let $X$ be an L-embedded Banach space with L-decomposition $X''=X\oplus_1 X_s$
and with L-projection $P$ from $X''$ onto $X$.
Let $(x_n)$ be a sequence spanning an almost isometric copy $Y$
of $\leins$ in $X$,
put $Y_m=[x_n]_{n\geq m}$, $m\in\N$.
Via the isomorphism between $Y$ and $\leins$
induced by $x_n\mapsto e_n$ the situation described for $\leins$
carries over to $Y$. That is, there is
$Z\subset Y^{\bot\bot}\subset X''$ such that
$Y^{\bot\bot}=Y\oplus Z$,
$Y_m^{\bot\bot}=Y_m\oplus (Z\cap Y_m^{\bot\bot})$ and
$z\in Z\cap Y_m^{\bot\bot}$ for all $m\in\N$, $z\in Z$.

Since $(x_n)$ is supposed to span $Y$ almost isometrically, there are
numbers $\eta_m\geq 0$ tending to $0$ such that
$\norm{y+z}\geq (1-\eta_m)(\norm{y}+\norm{z})$ for all $y\in Y_m$,
$z\in Z\cap Y_m^{\bot\bot}$.

Let $z\in Z$.
In order to show that $Y$ is L-embedded it is enough to show that
$Pz=0$ because then Li's criterion $PY^{\bot\bot}=Y$ is fulfilled.
But $z\in Z\cap Y_m^{\bot\bot}$ and
by a quantitative version of Li's result (see \cite[IV.1.4]{HWW}
or \cite[Lem.\ 2]{Pf-L1}) applied to $X$ and $Y_m$,
we have $\norm{Pz}\leq3\eta_m^{1/2}\norm{z}$
for all $m\in\N$. Hence $Pz=0$.

In passing we note (for Corollary \ref{coro_2} below) that
a \wst-accumulation point of $\{e_n\mdE n\in\N\}$ belongs to
$c_0^{\bot}\cap U_m^{\bot\bot}$ and has norm one. Accordingly,
if $z$ is a \wst-accumulation point of $\{x_n\mdE n\in\N\}$, then
$z\in X_s\cap Y_m^{\bot\bot}$ and $\norm{z}=1$,
the latter because each $Y_m^{\bot\bot}$ is $(1-\delta_m)$-isomorphic
to $(\leins)''$ with $\delta_m\gen0$.
\ebew\bigskip\\
In addition to Lemma \ref{prop_LLasy} we note that if $Y$ is an
almost isometric $\leins$-copy in an L-embedded space $X$ then
$X/Y$ is L-embedded, too, by \cite[IV.1.3]{HWW}.

Here is a way of constructing asymptotically isometric $\leins$-copies
in L-embedded Banach spaces. It is essentially due to Godefroy
\cite[IV.2.5]{HWW}:
\begin{satz}\label{lem_Go}
Let $X$ be an L-embedded Banach space with L-decomposition
$X''=X \oplus_1 X_s$,
$(\Gamma,\preceq)$ a directed set and
$(x_{\gamma})_{\gamma\in \Gamma}$
a net in the unit ball of $X$.
If $x_{\gamma}\genwst x_s\in X_s$ in the \wst-topology of $X''$
and $\norm{x_s}=1$
then there is a sequence
$(x_{\gamma_{n}})_{n\in N}$ which
spans $\leins$ asymptotically.
\end{satz}
{\em Proof}:
Let $(\delta_n)$ be a sequence of strictly positive numbers
converging to $0$.
Set $\eta_1=\frac{1}{3}\delta_1$ and
$\eta_{n+1}=\frac{1}{3}\min(\eta_n,\delta_{n+1})$ for
$n\in\N$.
By induction over $n\in\N$ we will construct $\gamma_n\in\Gamma$
such that 
\bgl
\Bigl(\sum_{i=1}^{n}(1-\delta_i)\betr{\alpha_i}\Bigr)
+\eta_n \sum_{i=1}^{n}\betr{\alpha_i}
\leq
\norm{\sum_{i=1}^{n}\alpha_i x_{\gamma_{i}}}
\leq         \sum_{i=1}^{n}\betr{\alpha_i}
\falle n\in\N, \,\,\alpha_i\in\C.
\label{glGo1}
\egl
(The last inequality is trivial because $\norm{x_{\gamma}}\leq 1$.)
We recall that the norm is \wst-lower semicontinuous.
Thus $\liminf\norm{x_{\gamma}}\geq \norm{x_s}=1$
and for the first induction step we
choose $x_{\gamma_{1}}$ such that
$\norm{x_{\gamma_{1}}}>1-\delta_1+\eta_1$.\\
For the induction step $n\mapsto n+1$ we suppose
$x_{\gamma_{1}},\ldots, x_{\gamma_{n}}$ to be constructed such that \Ref{glGo1}
holds.
Fix an element $\alpha=(\alpha_i)_{i=1}^{n+1}$ in the
unit sphere of $\leins_{n+1}$ such that $\alpha_{n+1}\neq 0$.
The \wst-convergence (along $\gamma$) of
$(\sum_{i=1}^n \alpha_i x_{\gamma_{i}}) + \alpha_{n+1}x_{\gamma}$ to
$(\sum_{i=1}^n \alpha_i x_{\gamma_{i}}) + \alpha_{n+1}x_s$ yields
\bglst
\liminf_{\gamma}
\Norm{\Bigl(\sum_{i=1}^n \alpha_i x_{\gamma_{i}}\Bigr) 
+ \alpha_{n+1}x_{\gamma}}
&\geq&
\Norm{\Bigl(\sum_{i=1}^n \alpha_i x_{\gamma_{i}}\Bigr) 
+ \alpha_{n+1}x_s}\\
&=&\Norm{\sum_{i=1}^n \alpha_i x_{\gamma_{i}}} + \norm{\alpha_{n+1}x_s}\\
&\stackrel{\Ref{glGo1}}{\geq}&
\Bigl(\sum_{i=1}^{n}(1-\delta_i)\betr{\alpha_i}\Bigr)
         +\eta_n \Bigl(\sum_{i=1}^{n}\betr{\alpha_i}\Bigr)
         +\betr{\alpha_{n+1}}\\
&=& \Bigl(\sum_{i=1}^{n+1}(1-\delta_i)\betr{\alpha_i}\Bigr)
        +\eta_n \Bigl(\sum_{i=1}^{n+1}\betr{\alpha_i}\Bigr)\\
&& \,\,\,\,\,\,\,\,\,\,\,\,\,\,\,\,\,\,\,\,\,\,\,\,
        -(\eta_n -\delta_{n+1})\betr{\alpha_{n+1}}\\
&\geq& \Bigl(\sum_{i=1}^{n+1}(1-\delta_i)\betr{\alpha_i}\Bigr)
        +\min(\eta_n,\delta_{n+1})
\eglst
because $\norm{\alpha}=1$ and $\betr{\alpha_{n+1}}\leq1$.
Thus there is an index $\gamma_0\in\Gamma$ such that
$\inf_{\gamma\succeq \gamma_0}
\Norm{\Bigl(\sum_{i=1}^n \alpha_i x_{\gamma_{i}}\Bigr) 
   + \alpha_{n+1}x_{\gamma}}
>    \Bigl(\sum_{i=1}^{n+1}(1-\delta_i)\betr{\alpha_i}\Bigr)+2\eta_{n+1}$;
note that the subnet $(x_{\gamma})_{\gamma\succeq \gamma_0}$
still \wst-converges to $x_s$.

Choose  a finite $\eta_{n+1}$-net
$(\alpha^{l})_{l=1}^{L_{n+1}}$ in the unit sphere
of $\leins_{n+1}$
in the sense that for each $\alpha$ in the unit sphere
of $\leins_{n+1}$ there is $l\leq L_{n+1}$ such that
$\norm{\alpha-\alpha^{l}} = \sum_{i=1}^{n+1}\betr{\alpha_i-\alpha_i^{l}}
<\eta_{n+1}$.
Then we may repeat the reasoning above finitely many times for
$l=1,\ldots, L_{n+1}$
in order to get $x_{\gamma_{n+1}}$ such that
\bglst
\Norm{\sum_{i=1}^{n+1}\alpha_i^{l} x_{\gamma_{i}}}
>
\Bigl(\sum_{i=1}^{n+1}(1-\delta_i)\betr{\alpha_i^{l}}\Bigr)+2\eta_{n+1}
\,\,\,\falle l\leq L_{n+1}.
\eglst
For an arbitrary $\alpha$ in the unit sphere
of $\leins_{n+1}$ choose $l\leq L_{n+1}$ such that
$\norm{\alpha-\alpha^{l}}<\eta_{n+1}$.
Then
\bgl
\Norm{\sum_{i=1}^{n+1}\alpha_i x_{\gamma_{i}}}
&\geq&
\Norm{\sum_{i=1}^{n+1}\alpha_i^{l} x_{\gamma_{i}}} -
   \Norm{\sum_{i=1}^{n+1}(\alpha_i - \alpha_i^{l}) x_{\gamma_{i}}}
\nonumber\\
&\geq&
\Bigl(\sum_{i=1}^{n+1}(1-\delta_i)\betr{\alpha_i}\Bigr)+2\eta_{n+1}
-\norm{\alpha-\alpha^{l}}
\nonumber\\
&\geq&
\Bigl(\sum_{i=1}^{n+1}(1-\delta_i)\betr{\alpha_i}\Bigr)+\eta_{n+1}
\nonumber\\
&=&
\Bigl(\sum_{i=1}^{n+1}(1-\delta_i)\betr{\alpha_i}\Bigr)
      +\eta_{n+1}\sum_{i=1}^{n+1}\betr{\alpha_i}.
\label{glGo2}
\egl
This extends to all $\alpha\in\leins_{n+1}$ and thus
ends the induction and the proof.
\ebew\bigskip\\
Within L-embedded spaces James' distortion theorem is more efficient:
\begin{coro}\label{coro_2}
A sequence spanning $\leins$ almost isometrically in an
L-embedded Banach space admits of a subsequence which spans
$\leins$ asymptotically.
\end{coro}
{\em Proof}:
We use the notation of the proof of Lemma \ref{prop_LLasy} and
the remark in the end of it.
Hence there is a net $(x_{n_{\gamma}})$ which \wst-converges to
a normalized \wst-accumulation point $x_s$ of $\{x_n\mdE n\in\N\}$.
Proposition \ref{lem_Go} applies.\ebew
\begin{coro}\label{coro_1}
Every nonreflexive subspace of an L-embedded Banach space
contains an asymptotic copy of $\leins$.\\
In particular, every nonreflexive subspace of an
L-embedded Banach space fails the fixed point property.
\end{coro}
{\em Proof}:
L-embedded spaces are \w-sequentially complete (\cite{G-bien}, or
\cite[IV.2.2]{HWW}). Hence the first assertion follows from
theorems of Eberlein-\v{S}muljan, Rosenthal, James and from
Corollary \ref{coro_2}.
Combine this with \cite[Th.\ 1.2]{DowLen} to get
the second assertion.
\ebew\bigskip\\
By \cite{Pf-L1} (or \cite[IV.2.7]{HWW})
the basis $(x_{\gamma_n})$ in Proposition \ref{lem_Go}
admits a subsequence whose closed linear span is complemented in $X$.
In this way Corollary \ref{coro_1} recovers
the theorem of Kadec and \Pel\
mentioned in the introduction.\medskip

From Lemma \ref{prop_LLasy} we see that both almost isometric and
asymptotically isometric $\leins$-copies are L-embedded
provided they are contained in an L-embedded Banach space.
But at the time of this writing it is not clear whether asymptotic
$\leins$-copies are always L-embedded (or whether asymptotic
$c_0$-copies are M-embedded).
For almost isometric  $\leins$-copies the situation is clearer:
\begin{coro}\label{coro_3}
There are almost isometric copies of $\leins$ which are not
L-embedded.
\end{coro}
{\em Proof}:
Combine Corollary \ref{coro_2} or \ref{coro_1}
and the counterexample of \cite{DJLT}.
\ebew\bigskip\\
Now we briefly turn to $c_0$.
Some straightforward
modifications of the proof of \cite{Pf-L1}
(or \cite[IV.2.7]{HWW}) show that the dual
of a nonreflexive L-embedded space contains
asymptotic copies of $c_0$ although, of course, not every
nonreflexive subspace in the dual of an L-embedded space
contains $c_0$-copies. (Take $l^{\infty}$ for example.)
M-embedded spaces provide a more natural frame:
\begin{proposition}\label{lem_c0}
Every nonreflexive subspace of an M-embedded space contains an
asymptotic $c_0$-copy.
\end{proposition}
{\em Proof}:
The result could be obtained by appropriate modifications
of the proof of \cite[III.3.4, 3.7a]{HWW} but we suggest
a more direct (and shorter) argument.

Since M-embeddedness passes to subspaces it is enough to prove
that every nonreflexive M-embedded space $Z$ contains an
asymptotic $c_0$-copy.
Let $(\delta_m)$ be a sequence in $]0,1[$ converging to $0$.
By induction over $n\in\N$ we will construct
a sequence $(z_n)$ in $Z$ such that
\bgl
\max_{i\leq n}(1-(1-2^{-n})\delta_i)\betr{\alpha_i}
\leq \norm{\sum_{i=1}^{n}\alpha_i z_i}
\leq \max_{i\leq n}(1+(1-2^{-n})\delta_i)\betr{\alpha_i}
\forall \alpha_i\in\C.
\label{gl_c01}
\egl
For the beginning of the induction
we choose $z_1$ in the unit sphere of $Z$.
Suppose that $z_1,\ldots, z_n$ have been constructed and satisfy
\Ref{gl_c01}. Let $P$ denote the L-projection on $Z'''$ with range $Z'$.
We have $P'\eing{Z}=\id_{Z}$ and
$Z''''=Z^{\bot\bot}\oplus_{\infty} Z'{}^{\bot}$.
There exists an element $z'{}^{\bot} \in Z'{}^{\bot}$
with $\norm{z'{}^{\bot}}=1$ because $Z$ is not reflexive.
Put $E=\mbox{lin}(\{z_i\mdE i\leq n\}\cup\{z'{}^{\bot}\})$ and choose
$\eta>0$ such that
\bglst
(1-\eta)(1-(1-2^{-n})\delta_i) &>& (1-(1-2^{-(n+1)})\delta_i)\\
(1+\eta)(1+(1-2^{-n})\delta_i) &<& (1+(1-2^{-(n+1)})\delta_i)
\eglst
for all $i\leq n$.
The principle of local reflexivity provides an operator
$T_1:E\gen Z''$ and an operator $T_2:T_1(E)\gen Z$
such that $T=T_2\circ T_1$ fulfills
\bgl
T\eing{E\cap Z}=\id_{E\cap Z}
\,\,\,\,\,\mbox{and}\,\,\,\,\,
(1-\eta)\norm{e}\leq\norm{Te}\leq(1+\eta)\norm{e}
\,\,\,\falle e\in E. \label{gl_c02}
\egl
Put $z_{n+1}=Tz'{}^{\bot}$. Then we get (\ref{gl_c01}, $n+1$):
\bglst
\max_{i\leq n+1}(1-(1-2^{-(n+1)})\delta_i)\betr{\alpha_i}
&<&
(1-\eta)\max\Bigl(\max_{i\leq n}(1-(1-2^{-n})\delta_i)\betr{\alpha_i},
              \betr{\alpha_{n+1}}\Bigr)\\
&\stackrel{\Ref{gl_c01}}{\leq}&
(1-\eta)\max\Bigl(\norm{\sum_{i=1}^{n}\alpha_i z_i},
       \norm{\alpha_{n+1}z'{}^{\bot}}\Bigr)\\
&=&(1-\eta)\norm{\Bigl(\sum_{i=1}^{n}\alpha_i z_i\Bigr)
                          +\alpha_{n+1}z'{}^{\bot}}\\
&\stackrel{\Ref{gl_c02}}{\leq}&
       \norm{\sum_{i=1}^{n+1}\alpha_i z_i}\\
&\stackrel{\Ref{gl_c02}}{\leq}&
(1+\eta)\norm{\Bigl(\sum_{i=1}^{n}\alpha_i z_i\Bigr)
                          +\alpha_{n+1}z'{}^{\bot}}\\
&=&(1+\eta)\max\Bigl(\norm{\sum_{i=1}^{n}\alpha_i z_i},
       \betr{\alpha_{n+1}}\Bigr)\\
&<&\max_{i\leq n+1}(1+(1-2^{-(n+1)})\delta_i)\betr{\alpha_i}.
\eglst
This ends the induction and the proof.
\ebew\medskip

It is now immediate that every nonreflexive subspace of
an M-embedded space fails the fixed point property
(see \cite[Prop.\ 7]{DowLenTur}).

The analogue of Lemma \ref{prop_LLasy} holds trivially,
simply because M-embeddedness passes to subspaces \cite[III.1.6]{HWW};
thus every $c_0$-copy inside an M-embedded space is M-embedded.
And similarly as for Corollary \ref{coro_3}, 
the counterexample of \cite{DJLT}
combined with Proposition \ref{lem_c0} provides an
almost isometric  $c_0$-copy which is not M-embedded.

\end{document}